\documentclass{amsart}
\usepackage{graphics}

\newtheorem{theorem}{Theorem}[section]
\newtheorem{lemma}[theorem]{Lemma}

\theoremstyle{definition}

\theoremstyle{remark}

\numberwithin{equation}{section}

\DeclareSymbolFont{AMSb}{U}{msb}{m}{n}
\DeclareMathSymbol{\Z}{\mathalpha}{AMSb}{"5A}

\begin{document}
\newcommand{\beqs}{\begin{equation*}}
\newcommand{\eeqs}{\end{equation*}}
\newcommand{\beq}{\begin{equation}}
\newcommand{\eeq}{\end{equation}}
\newcommand\nutwid{\overset {\text{\lower 3pt\hbox{$\sim$}}}\nu}
\newcommand\Mtwid{\overset {\text{\lower 3pt\hbox{$\sim$}}}M}
\newcommand\ptwid{\overset {\text{\lower 3pt\hbox{$\sim$}}}p}
\newcommand\pitwid{\overset {\text{\lower 3pt\hbox{$\sim$}}}\pi}
\newcommand\bijone{\overset {1}\longrightarrow}                   
\newcommand\bijtwo{\overset {2}\longrightarrow}                   
\newcommand\pihat{\widehat{\pi}}
\newcommand\mymod[1]{(\mbox{mod}\ {#1})}
\newcommand\myto{\to}
\newcommand\srank{\mathrm{srank}}
\newcommand\pitc{\pi_{\mbox{$t$-core}}}
\newcommand\pitcb[1]{\pi_{\mbox{${#1}$-core}}}  
\newcommand\mpitcb[1]{\pi_{\mathrm{{#1}-core}}}  
\newcommand\lamseq{(\lambda_1, \lambda_2, \dots,\lambda_\nu)}
\newcommand\mylabel[1]{\label{#1}}
\newcommand\eqn[1]{(\ref{eq:#1})}
\newcommand\stc{{St-crank}}
\newcommand\mstc{\mbox{St-crank}}
\newcommand\mmstc{\mbox{\scriptsize\rm St-crank}}
\newcommand\tqr{{$2$-quotient-rank}}
\newcommand\mtqr{\mbox{$2$-quotient-rank}}
\newcommand\mmtqr{\mathrm{2-quotient-rank}}
\newcommand\fcc{{$5$-core-crank}} 
\newcommand\fc{{5\mbox{-core}}} 
\newcommand\tc{{t\mbox{-core}}}
\newcommand\pifc{\pi_{5\mbox{-core}}} 
\newcommand\bgr{BG-rank} 
\newcommand\mbgr{\mbox{BG-rank}} 
\newcommand\nvec{(n_0, n_1, \dots, n_{t-1})}
\newcommand\parity{\mbox{par}}
\newcommand\epart{\mbox{ep}}
\newcommand\pbar{\overline{p}}
\newcommand\Pbar{\overline{P}}
\newcommand\abar{\overline{a}}
\newcommand\legendre[2]{\genfrac{(}{)}{}{}{#1}{#2}}

\title[The BG-rank of a partition and its applications]{The BG-rank of a partition
and its applications}

\author{Alexander Berkovich}
\address{Department of Mathematics, University of Florida, Gainesville,
Florida 32611-8105}
\email{alexb@math.ufl.edu}          
\thanks{Research of both authors was supported in part by NSA grant MSPF-06G-150.}

\author{Frank G. Garvan}
\address{Department of Mathematics, University of Florida, Gainesville,
Florida 32611-8105}
\email{frank@math.ufl.edu}          


\date{April 27, 2007}  



\begin{abstract}
Let $\pi$ denote a partition into parts $\lambda_1\ge\lambda_2\ge\lambda_3\ldots$.
In a 2006 paper 
we defined \bgr $(\pi)$ as 
\beqs
\mbgr(\pi) = \sum_{j\ge 1}(-1)^{j+1}\frac{1-(-1)^{\lambda_j}}{2}.
\eeqs
This statistic was employed to generalize and refine the famous Ramanujan modulo $5$ partition congruence.
Let $p_j(n)$
denote the number of partitions of $n$ with $\mbgr=j$. Here, we provide
a combinatorial proof that
\beqs
p_j(5n+4)\equiv 0 \pmod{5},\quad j\in \mathbb Z,
\eeqs
by showing that the residue of the $5$-core crank mod $5$ 
divides the partitions enumerated by $p_j(5n+4)$ into five equal classes. 
This proof uses the orbit construction from our previous paper and 
a new identity for the \bgr.
Let $a_{t,j}(n)$ denote the number of $t$-cores of $n$
with $\mbgr=j$. We find
eta-quotient representations for
\beqs
\sum_{n\ge 0}a_{t,\big\lfloor\frac{t+1}{4}\big\rfloor}(n)q^n \quad \mbox{and} \quad
\sum_{n\ge 0} a_{t,-\big\lfloor\frac{t-1}{4}\big\rfloor}(n)q^n,
\eeqs
when $t$ is an odd, positive integer. Finally, we derive explicit formulas for the coefficients $a_{5,j}(n)$, $j=0,\pm 1$.
\end{abstract}
\maketitle

\section{Introduction} \label{sec:intro}

A partition $\pi$ is a nonincreasing sequence
\beqs
\pi = (\lambda_1,\lambda_2,\lambda_3,\ldots)
\eeqs
of positive integers (parts) $\lambda_1\geq\lambda_2\geq\lambda_3\geq\ldots$.
The norm of $\pi$, denoted $|\pi|$, is defined as
\beqs
|\pi| = \sum_{i\geq 1}\lambda_i.
\eeqs
If $|\pi|=n$, we say that $\pi$ is a partition of $n$. The (Young) diagram of $\pi$
is a convenient way to represent $\pi$ graphically: the parts of $\pi$ are shown as rows
of unit squares (cells). Given the diagram of $\pi$ we label a cell in the $i$-th row
and $j$-th column by the least nonnegative integer $\equiv j-i\pmod{t}$. The resulting
diagram is called a $t$-residue diagram \cite{K}. We can also label cells in the infinite column 
$0$ and the infinite row $0$ in the same fashion and call the resulting diagram the 
extended $t$-residue diagram \cite{GKS}.
And so with each partition $\pi$ and positive integer $t$ we can associate the $t$-dimensional vector
\beqs
\vec{r}(\pi,t)=(r_0(\pi,t),r_1(\pi,t),\ldots,r_{t-1}(\pi,t))
\eeqs
with
\beqs
r_i(\pi,t)=r_i, \quad 0\leq i \leq t-1
\eeqs
being the number of cells colored $i$ in the $t$-residue diagram of $\pi$.
If some cell of $\pi$ shares a vertex or edge with the rim of the diagram of $\pi$, 
we call this cell a rim cell of $\pi$. A connected collection of rim cells of $\pi$ is
called a rim hook if (diagram of $\pi$)$\backslash$(rim hook) represents a legitimate partition. 
We say that a partition is a $t$-core, denoted $\pi_{t\mbox{-core}}$, 
if its diagram has no rim hooks of length $t$ \cite{K}.

The Durfee square of $\pi$ is the largest square that fits inside the diagram of $\pi$.
Reflecting the diagram of $\pi$ about its main diagonal, one gets the diagram of $\pi'$
(the conjugate of $\pi$). More formally, 
\beqs
\pi' = (\lambda_1',\lambda_2',\lambda_3',\ldots)
\eeqs
with $\lambda_i'$ being the number of parts of $\pi$ that are $\geq i$.
In \cite{BG} we defined a new partition statistic
\beq
\mbgr(\pi): = \sum_{j\ge 1}(-1)^j\frac{(-1)^{\lambda_j}-1}{2}.
\mylabel{eq:bgd} 
\eeq
It is easy to verify that
\beq
\mbgr(\pi)=r_0(\pi,2)-r_1(\pi,2)
\mylabel{eq:bgd'} 
\eeq
and
\beq
\mbgr(\pi)\equiv \vert\pi\vert\pmod 2.
\mylabel{eq:aform} 
\eeq
In \cite{BG} we proved the following $\pmod{5}$ congruences
\begin{align}
p_j(5n)   &\equiv 0 \pmod{5} \quad \mbox{if } j    \equiv 1,2\pmod{5},\mylabel{eq:10} \\ 
p_j(5n+1) &\equiv 0 \pmod{5} \quad \mbox{if } j\not\equiv 1,2\pmod{5},\mylabel{eq:20} \\ 
p_j(5n+2) &\equiv 0 \pmod{5} \quad \mbox{if } j\not\equiv 0,3\pmod{5},\mylabel{eq:30} \\ 
p_j(5n+3) &\equiv 0 \pmod{5} \quad \mbox{if } j    \equiv 0,3\pmod{5},\mylabel{eq:40} \\ 
p_j(5n+4) &\equiv 0 \pmod{5} \quad \mbox{for all } j\in\mathbb Z. \mylabel{eq:*} 
\end{align}
Here $p_j(n)$ denotes the number of partitions of $n$ with $\mbgr=j$.
Clearly,
\beqs
p(5n+4) = \sum_j p_j(5n+4)
\eeqs
with $p(n)$ denoting the number of unrestricted partitions of $n$. And so (\ref{eq:*})
implies the famous Ramanujan congruence \cite{R}
\beqs
p(5n+4)\equiv 0\pmod{5}.
\eeqs
In this paper, we build on the developments in \cite{BG} to provide a 
combinatorial proof of (\ref{eq:*}). 

For $t$-odd it is surprising that the $\mbgr(\pi_{t-\mbox{core}})$ assumes only finitely many values. In fact,
we will show that if $t$ is an odd, positive integer, then
\beq
-\bigg\lfloor\frac{t-1}{4}\bigg\rfloor\leq\mbgr(\pi_{t-\mbox{core}})\leq\bigg\lfloor\frac{t+1}{4}\bigg\rfloor.
\mylabel{eq:bound} 
\eeq
Here $\lfloor x\rfloor$ denotes the integer part of $x$. 

We will establish the following 
identities. For odd $t>1$
\begin{align}
C_{t,(-1)^{\frac{t-1}{2}}\big\lfloor\frac{t-1}{4}\big\rfloor}(q)  & =q^{\frac{(t-1)(t-3)}{8}}\quad F(t,q^2),
\mylabel{eq:**}\\ 
C_{t,(-1)^{\frac{t+1}{2}}\big\lfloor\frac{t+1}{4}\big\rfloor}(q)& =q^{\frac{t^2-1}{8}}\quad \frac{E^t(q^{4t})}{E(q^4)}, \mylabel{eq:***} 
\end{align}
where
\beqs
C_{t,j}(q)=\sum_{n\geq0}a_{t,j}(n)q^n,
\eeqs
$a_{t,j}(n)$ denotes the number of $t$-cores of $n$ with $\mbgr=j$ and
\beqs
E(q)=\prod_{j=1}^\infty (1-q^j),
\eeqs
\beqs
F(t,q)=\frac{E^{t-4}(q^{2t}) E^2(q^t) E^3(q^2)}{E^2(q)}.
\eeqs
We observe that (\ref{eq:aform}) suggests that $C_{t,j}(q)$ is an even (odd) function of $q$ if $j$ is even (odd).

It is instructive to compare (\ref{eq:**}, \ref{eq:***}) with the well-known identity \cite{GKS} 
for unrestricted $t$-cores
\beq
\sum_{n\ge0}a_t(n)q^n=\frac{E^t(q^t)}{E(q)}.
\mylabel{eq:+} 
\eeq
Here $a_t(n)$ denotes the number of $t$-cores of $n$.

The rest of this paper is organised as follows.

In Section 2 we discuss the Littlewood decomposition of $\pi$ in terms of $t$-core and $t$-quotient of $\pi$.
We describe the Garvan, Kim, Stanton bijection for $t$-cores and use a constant 
term technique to provide a simple proof of the Klyachko identity \cite{K1}
\beq
\sum_{\substack{\vec n\in \mathbb Z^t \\ \vec n\cdot\vec{1}_t=0}} q^{\frac{t}{2}\vec n\cdot\vec{n}+\vec b_t\cdot\vec{n}}=
\frac{E^t(q^t)}{E(q)}.
\mylabel{eq:++} 
\eeq
Here $\vec{1}_t=(1,1,\ldots,1)\in \mathbb Z^t$, $\vec b_t=(0,1,2,\ldots,t-1)$.

In Section 3 we establish a fundamental identity connecting $\mbgr$ and the Littlewood decomposition.

In Section 4 we discuss a combinatorial proof of (\ref{eq:*}).

Section 5 is devoted to the proof of the identities (\ref{eq:**}, \ref{eq:***}).

Section 6 deals with $5$-cores with prescribed $\mbgr$. There we derive the explicit formulas for the coefficients $a_{5,j}(n)$, $j=0,\pm 1$. 

In Section 7 we give a generalization of the \bgr\ and state a number of results.

\section{Two Bijections} \label{sec:2bj}

In this section we will follow closely the discussion in  \cite{G}, \cite{GKS} to recall some basic facts 
about $t$-cores and $t$-quotients. A region $r$ in the extended $t$-residue diagram of $\pi$ is the set 
of cells $(i,j)$ satisfying $t(r-1)\leq j-i<tr$. A cell of $\pi$ is called exposed if it is at the end of a row. 
One can construct $t$ bi-infinite words $W_0,W_1,\ldots,W_{t-1}$ of two letters $N,E$ as
\begin{align*}
\mbox{The $r$-th letter of }W_i=
\begin{cases}
E, &\mbox{if there is an exposed cell labelled $i$ in the region $r$}\\
N, &\mbox{otherwise}.
\end{cases}
\end{align*}
It is easy to see that the word set $\{W_0,W_1,\ldots,W_{t-1}\}$ fixes $\pi$ uniquely.

Let $P$ be the set of all partitions and $P_{t\mbox{-core}}$ be the set of all $t$-cores.
There is a well-known bijection
\beqs
\phi_1:P\rightarrow P_{t\mbox{-core}}\times P \times P \times P \ldots\times P 
\eeqs
which goes back to Littlewood \cite{L}
\beqs
\phi_1(\pi)=(\pi_{t\mbox{-core}},\hat\pi_0,\hat\pi_1,\ldots,\hat\pi_{t-1})
\eeqs
such that
\beqs
|\pi|=|\pi_{t\mbox{-core}}|+t\sum_{i=0}^{t-1}|\hat\pi_i|.
\eeqs
Multipartition $(\hat\pi_0,\hat\pi_1,\ldots,\hat\pi_{t-1})$ is called the $t$-quotient of $\pi$. 
We remark that (\ref{eq:+}) is the immediate corollary of the Littlewood bijection. We describe $\phi_1$
in full detail a bit later.

The second bijection 
\beqs
\phi_2:P_{t\mbox{-core}}\rightarrow \{\vec{n}:\vec{n}\in\mathbb Z^t,\vec{n}\cdot\vec 1_t=0\}
\eeqs
was introduced in \cite{GKS}. It is for $t$-cores only
\beqs
\phi_2(\pi_{t\mbox{-core}})=\vec{n}=(n_0,n_1,\ldots,n_{t-1})
\eeqs
where for $0\leq i\leq t-2$
\beq
n_i=r_i(\pi_{t\mbox{-core}},t)-r_{i+1}(\pi_{t\mbox{-core}},t)
\mylabel{eq:nr} 
\eeq
and
\beq
n_{t-1}=r_{t-1}(\pi_{t\mbox{-core}},t)-r_0(\pi_{t\mbox{-core}},t).
\mylabel{eq:nr'} 
\eeq
Clearly, 
\beqs
\sum_{i=0}^{t-1}n_i=\vec n \cdot \vec 1_t=0.
\eeqs
Moreover,
\beq
|\pi_{t\mbox{-core}}|=\frac{t}{2}\vec{n}\cdot\vec{n}+\vec{b}_t\cdot\vec n,
\mylabel{eq:R} 
\eeq
as shown in \cite{GKS}. And so
\beq
\sum_{n\geq 0}a_t(n)q^n=
\sum_{\substack{\vec n\in \mathbb Z^t \\ \vec n\cdot\vec{1}_t=0}} 
q^{\frac{t}{2}\vec n\cdot\vec{n}+\vec b_t\cdot\vec{n}}.
\mylabel{eq:X} 
\eeq
Note that (\ref{eq:+}), (\ref{eq:X}) imply the Klyachko identity (\ref{eq:++}).
The reader may wonder if (\ref{eq:nr}, \ref{eq:nr'}) can be used to define $\phi_2(\pi)=\vec n$ 
for any partition $\pi$. This, of course, can be done. However, in general $\phi_2$ is not a $1-1$ 
function and so $\phi_2^{-1}$ can't be defined. Indeed, if $\pi_1\neq\pi_2$, but $\pi_{t\mbox{-core}}$
is a $t$-core of both $\pi_1$ and $\pi_2$ then
\beqs
\phi_2(\pi_1)=\phi_2(\pi_2)=\phi_2(\pi_{t\mbox{-core}}).
\eeqs
When a partition is a $t$-core, $\phi_2$ can be inverted. To do this we recall that the partition is a $t$-core iff for 
$0\leq i\leq t-1$
\beqs
\begin{matrix}
\mbox{Region}&: &\cdots\cdots\cdots & n_{i}-1 & n_i & n_{i}+1 & n_{i}+2 & \cdots\cdots\cdots \\
          W_i&: &\cdots\cdots\cdots & E       & E   & N       & N       & \cdots\cdots\cdots 
\end{matrix}
\eeqs
as explained in \cite{GKS}. For example, the word image of $\phi_2^{-1}((2,-1,-1))$ is
\begin{align*}
\mbox{Region}&:\quad\cdots\cdots -1 \quad\; 0 \quad\; 1 \quad\; 2 \quad\; 3 \quad\cdots\cdots \\
W_0&:\quad\cdots\cdots\quad E \quad E \quad E \quad E \quad N \quad\cdots\cdots \\
W_1&:\quad\cdots\cdots\quad E \quad N \quad N \quad N \quad N \quad\cdots\cdots \\
W_2&:\quad\cdots\cdots\quad E \quad N \quad N \quad N \quad N \quad\cdots\cdots.
\end{align*}
This means that
\beq
\phi_2^{-1}((2,-1,-1))=(4,2).
\mylabel{eq:B} 
\eeq
More generally, if
\beqs
\phi_1(\pi)=(\pi_{t\mbox{-core}},\hat\pi_0,\hat\pi_1,\ldots,\hat\pi_{t-1})
\eeqs
with
\beqs
\hat\pi_i=(\lambda_1^{(i)},\lambda_2^{(i)},\ldots,\lambda_{m_i}^{(i)}), \quad 0\leq i \leq t-1,
\eeqs
then cells colored $i$ are not exposed only in the regions
\beqs
n_i+j-\lambda_j^{(i)}, \quad 1\leq j\leq m_i
\eeqs
and
\beqs
n_i+m_i+k, \quad k\geq 1.
\eeqs
For example, if $\hat\pi_i=(\lambda_1,\lambda_2,\lambda_3)$ then
\begin{align*}
\mbox{Region}&:\quad\cdots\cdots &n_i+1-\lambda_1 \quad\cdots\cdots\quad &n_i+2-\lambda_2
               \quad\cdots\cdots &n_i+3-\lambda_3 \quad\cdots\cdots\quad &n_i+4 \quad\cdots\cdots\\ 
          W_i&:\quad\cdots\cdots &E \quad N \quad E \quad\cdots\cdots\quad &E \quad N \quad E 
               \quad\cdots\cdots &E \quad N \quad E \quad\cdots\cdots\quad &E \quad N \quad\cdots\cdots
\end{align*}
Clearly, one can easily determine $\vec n$ and $(\hat\pi_0,\hat\pi_1,\ldots,\hat\pi_{t-1})$ from the word
set $\{W_0,W_1,\ldots,W_{t-1}\}$. And so
\beqs
\phi_1(\pi)=(\phi_2^{-1}(\vec n),\hat\pi_0,\ldots,\hat\pi_{t-1}).
\eeqs
We illustrate the above with the following example. If $t=3$ and $\pi=(7,5,4,3,2)$ then
\begin{align*}
\mbox{Region}&:\quad\cdots\cdots\; -2 \; -1\quad 0\quad 1\quad\;\; 2\quad\;\; 3\quad\; 4\quad\; 5\quad\cdots\cdots \\
          W_0&:\quad\cdots\cdots\quad E \quad E \quad E \quad N \quad E \quad E \quad N \quad N \quad\cdots\cdots \\
          W_1&:\quad\cdots\cdots\quad E \quad N \quad N \quad E \quad N \quad N \quad N \quad N \quad\cdots\cdots \\
          W_2&:\quad\cdots\cdots\quad E \quad N \quad E \quad N \quad N \quad N \quad N \quad N \quad\cdots\cdots.
\end{align*}
We have
\begin{align*}
n_0&=2,  &\pihat_0&=(2),\\
n_1&=-1, &\pihat_1&=(1,1),\\
n_2&=-1, &\pihat_2&=(1).
\end{align*}
Using (\ref{eq:B}), we obtain
\beqs
\phi_1((7,5,4,3,2))=((4,2),(2),(1,1),(1)).
\eeqs
To proceed further we recall some standard $q$-hypergeometric notations \cite{GR}:
\beqs
(a_1,a_2,a_3,\ldots ;q)_N=(a_1;q)_N (a_2;q)_N (a_3;q)_N \ldots
\eeqs
where
\begin{align*}
(a;q)_N=(a)_N=
\begin{cases}
\prod_{j=0}^{N-1}(1-aq^j), & N>0 \\
1,                         & N=0 \\
\prod_{j=1}^{-N}(1-aq^{-j})^{-1}, & N<0.
\end{cases}
\end{align*}
We shall also require the Jacobi triple product identity \cite[(II.28)]{GR}
\beq
\sum_{n=-\infty}^\infty q^{n^2} z^n = (q^2,-zq,-\frac{q}{z};q^2)_\infty.
\mylabel{eq:J} 
\eeq
We are now ready to prove the Klyachko identity (\ref{eq:++}). We will employ a so-called constant term technique.
To this end we rewrite the left hand side of (\ref{eq:++}) as
\beqs
\mbox{LHS (\ref{eq:++}) }=[z^0] \sum_{\vec n \in\mathbb Z^t} 
q^{\frac{t}{2}\vec n\cdot\vec n+\vec b_t\cdot\vec n} z^{\vec n\cdot\vec 1_t}=
[z^0] \prod_{i=0}^{t-1}\sum_{n_i=-\infty}^\infty 
q^{\frac{t}{2} n_i^2 + in_i} z^{n_i}
\eeqs
where $[z^i]f(z)$ is the coefficient of $z^i$ in the expansion of $f(z)$ in powers of $z$.
With the aid of (\ref{eq:J}) we derive
\begin{align*}
\mbox{LHS (\ref{eq:++}) } 
&=[z^0]\prod_{i=0}^{t-1}\left(q^t,-q^{i+\frac{t}{2}}z,-\frac{q^{\frac{t}{2}}}{q^iz};q^t\right)_\infty \\
&=[z^0]\frac{E^t(q^t)}{E(q)}\left(q,-q^{\frac{t}{2}}z,-\frac{q}{q^{\frac{t}{2}}z};q\right)_\infty \\
&=[z^0]\bigg(\frac{E^t(q^t)}{E(q)}\sum_{n=-\infty}^\infty q^{\frac{n^2}{2}+\frac{t-1}{2}n}z^n\bigg) \\
&=\frac{E^t(q^t)}{E(q)},
\end{align*}
as desired. The above proof is just a warm-up excercise to prepare the reader for a more sophisticated proof 
of (\ref{eq:**}) discussed in Section 5.

\section{The Littlewood decomposition and \bgr} \label{sec:3}

The main goal of this section is to establish the following 
identities for $\mbgr$.
If $t$ is even and $(n_0,\ldots,n_{t-1})=\phi_2(\pi)$, then 
\beq
\mbgr(\pi)=\sum_{i=0}^{\frac{t-2}{2}}n_{2i}.
\mylabel{eq:1} 
\eeq
If $t$ is odd then
\beq
\mbgr(\pi_{t\mbox{-core}})=bg(\vec n),
\mylabel{eq:2} 
\eeq
where $\vec n=\phi_2(\pi_{t\mbox{-core}})$ and
\beq
bg(\vec n):=\frac{1-\sum_{j=0}^{t-1}(-1)^{j+n_j}}{4}.
\mylabel{eq:3} 
\eeq
Moreover, if $t$ is odd and $\phi_1(\pi)=(\pi_{t\mbox{-core}},\pihat_0,\ldots,\pihat_{t-1})$ then
\beq
\mbgr(\pi)=\mbgr(\pi_{t\mbox{-core}})+\sum_{j=0}^{t-1}(-1)^{j+n_j}\mbgr(\pihat_j).
\mylabel{eq:4} 
\eeq
The proof of (\ref{eq:1}) is straightforward. It is sufficient to observe that if some cell is colored $i$ in the $t$-residue diagram of $\pi$, then it is colored $\frac{1-(-1)^i}{2}$ in the $2$-residue diagram of $\pi$. And so we obtain with the aid of (\ref{eq:bgd'})
\begin{align*}
\mbgr(\pi)&=(r_0+r_2+r_4+\cdots+r_{t-2})-(r_1+r_3+r_5+\cdots+r_{t-1}) \\
          &=(r_0-r_1)+(r_2-r_3)+\cdots+(r_{t-2}-r_{t-1}) \\
          &=n_0+n_2+\cdots+n_{t-2},
\end{align*}
as desired. Next, let $D(\pi)=D$ denote the size of the Durfee square of $\pi$. To prove (\ref{eq:2}) we begin by rewriting (\ref{eq:bgd}) as
\beq
\mbgr(\pi)=\frac{1}{2}\left(\parity(\nu)+\sum_{j=1}^\nu(-1)^{\lambda_j-j}\right).
\mylabel{eq:5} 
\eeq
Here $\pi=(\lambda_1,\lambda_2,\ldots,\lambda_\nu)$ and $\parity(x)$ is defined as
\beqs
\parity(x):=\frac{1-(-1)^x}{2}.
\eeqs
Next, let $\pi_1,\pi_2$ denote the partitions constructed from the first $D=D(\pi_{t\mbox{-core}})$ rows,
columns of $\pi_{t\mbox{-core}}$, respectively. Let $\pi_3$ denote a partition whose diagram is the Durfee square of
$\pi_{t\mbox{-core}}$. It is plain that
\begin{align}
\mbgr(\pi_{t\mbox{-core}})&=\mbgr(\pi_1)+\mbgr(\pi_2)-\mbgr(\pi_3)\nonumber \\
&=\mbgr(\pi_1)+\mbgr(\pi_2)-\parity(D).
\mylabel{eq:6} 
\end{align}
We shall also require the following sets
\begin{align*}
P_+:&= \{i\in \mathbb Z: \quad 0\leq i\leq t-1, \quad n_i>0 \}, \\
P_-:&= \{i\in \mathbb Z: \quad 0\leq i\leq t-1, \quad n_i<0 \}.
\end{align*}
Here $n_i$'s are the components of $\phi_2(\pi_{t\mbox{-core}})$. Note that if $i\in P_+$, then $i$ is exposed in all positive regions $\leq n_i$ of $\pi_1$. This observation together with (\ref{eq:5}) implies that
\begin{align}
\mbgr(\pi_1)&=\frac{1}{2}\left(\parity(D)+\sum_{i\in P_+}\sum_{k=1}^{n_i}(-1)^{t(k-1)+i}\right) \nonumber\\
						&=\frac{1}{2}\left(\parity(D)+\sum_{i\in P_+}(-1)^i \parity(n_i)\right)
\mylabel{eq:8} 
\end{align}

In \cite{GKS}, the authors showed that under conjugation $\phi_2(\pi_{t\mbox{-core}})$ transforms as
\beqs
(n_0,n_1,n_2,\ldots,n_{t-1})\rightarrow(-n_{t-1},-n_{t-2},-n_{t-3},\ldots,-n_0).
\eeqs
Also it is easy to see that
\beqs
\mbgr(\pi_2)=\mbgr(\pi_2').
\eeqs
It follows that
\beq
\mbgr(\pi_2)=\frac{1}{2}\left(\parity(D)+\sum_{i\in P_-}(-1)^i \parity(n_i)\right).
\mylabel{eq:9} 
\eeq
Combining (\ref{eq:6}, \ref{eq:8}, \ref{eq:9}) and taking into account that $\parity(0)=0$ we get
\begin{align*}
\mbgr(\pi_{t\mbox{-core}})&=\frac{1}{2}\sum_{i\in P_-\cup P_+}(-1)^i \parity(n_i) \\
                          &=\frac{1}{2}\sum_{i=0}^{t-1}(-1)^i \parity(n_i)=\frac{1-\sum_{i=0}^{t-1}(-1)^{i+n_i}}{4},
\end{align*}
as desired.
Note that formula (\ref{eq:2}) implies that $\mbgr$ of odd $t$-core is bounded, as stated in (\ref{eq:bound}).
Next, let $\pitwid_{0,i}, \pitwid_{2,i}, \pitwid_{3,i},\ldots$ denote the partitions constructed from 
$\phi_1(\pi)=(\pi_{t\mbox{-core}},\pihat_0,\pihat_1,\ldots,\pihat_{t-1})$, 
for odd $t$ as follows
\begin{align*}
\pitwid_{0,i}&=
\phi_1^{-1}(\pi_{t\mbox{-core}},\pihat_0,\pihat_1,\ldots,\pihat_{i-1},(0),\pihat_{i+1},\ldots,\pihat_{t-1}),\\
\pitwid_{1,i}&=
\phi_1^{-1}(\pi_{t\mbox{-core}},\pihat_0,\pihat_1,\ldots,\pihat_{i-1},(\lambda_1),\pihat_{i+1},\ldots,\pihat_{t-1}),\\
\pitwid_{2,i}&=
\phi_1^{-1}(\pi_{t\mbox{-core}},\pihat_0,\pihat_1,\ldots,\pihat_{i-1},(\lambda_1,\lambda_2),
                               \pihat_{i+1},\ldots,\pihat_{t-1}), \\
&\cdots\cdots\cdots.                               
\end{align*}
Here $\pihat_i=(\lambda_1,\lambda_2,\ldots,\lambda_\nu)$. Note that the $W_i$ word of $\pitwid_{0,i}$ is
\beqs
\begin{matrix}
\mbox{Region}&: &\cdots\cdots\cdots & n_i  & n_{i}+1 & \cdots\cdots\cdots \\
          W_i&: &\cdots\cdots\cdots & E    & N       & \cdots\cdots\cdots.
\end{matrix}
\eeqs
To convert $\pitwid_{0,i}$ into $\pitwid_{1,i}$ we attach a rim hook of length $t\lambda_1$ to $\pitwid_{0,i}$
so that $W_i$ becomes
\beqs
\begin{matrix}
\mbox{Region}&: &\cdots\cdots\cdots & n_i+1-\lambda_1   &\cdots\cdots\cdots & n_{i}+2,           & \cdots\cdots\cdots \\
          W_i&: &\cdots\cdots\cdots & E \quad N \quad E &\cdots\cdots\cdots & E \quad N \quad N  & \cdots\cdots\cdots.
\end{matrix}
\eeqs
It is not hard to verify that the color of the head (north-eastern) cell of the added rim-hook in the $2$-residue diagram of $\pitwid_{1,i}$ is given by $\parity(tn_i+i)=\parity(n_i+i)$. Observe that zeros and ones alternate along the added hook rim. This means that $\mbgr$ does not change if $\lambda_1$ is even. If $\lambda_1$ is odd then the change is determined by the color of the added head cell, i.e.
\begin{align*}
\mbgr(\pitwid_{1,i})&=\mbgr(\pitwid_{0,i})+\parity(\lambda_1)(1-2\parity(n_i+i))\\
                    &=\mbgr(\pitwid_{0,i})+\parity(\lambda_1)(-1)^{n_i+i},
\end{align*}
Next, we convert $\pitwid_{1,i}$ into $\pitwid_{2,i}$ by adding the new hook rim of length $t\lambda_2$ to $\pitwid_{1,i}$ so that $W_i$ becomes
\beqs
\begin{matrix}
\mbox{Region}&: &\cdots\cdots &n_i+1-\lambda_1 &\cdots\cdots & n_i+2-\lambda_2 &\cdots\cdots & n_i+3   &\cdots\cdots \\
          W_i&: &\cdots\cdots &E\quad N\quad E &\cdots\cdots & E\quad N\quad E &\cdots\cdots & E\quad N  &\cdots\cdots.
\end{matrix}
\eeqs
The color of the new head cell is given by
\beqs
\parity(t(n_i+1)+i)=\parity(n_i+1+i),
\eeqs
and so
\begin{align*}
\mbgr(\pitwid_{2,i})&=\mbgr(\pitwid_{1,i})+\parity(\lambda_2)(1-2\parity(n_i+1+i))\\
                    &=\mbgr(\pitwid_{0,i})+(-1)^{n_i+i}(\parity(\lambda_1)-\parity(\lambda_2)).
\end{align*}
Proceeding as above we arrive at
\begin{align}
\mbgr(\pi)&=\mbgr(\pitwid_{0,i})+(-1)^{n_i+i}\sum_{j=1}^\nu (-1)^{j+1}\parity(\lambda_j) \nonumber \\
          &=\mbgr(\pitwid_{0,i})+(-1)^{n_i+i}\mbgr(\pihat_i).
\mylabel{eq:310} 
\end{align}
Formula (\ref{eq:4}) follows easily from (\ref{eq:310}). Let us now define $\vec B_t,\vec{\tilde B}_t\in\mathbb Z^t$ as
\begin{align*}
\vec{B}_t=
\begin{cases}
\sum_{i=0}^{\frac{t-1}{2}} \vec e_{2i},   & \mbox{if } t\equiv  1\pmod 4 \\
\sum_{i=0}^{\frac{t-3}{2}} \vec e_{1+2i}, & \mbox{if } t\equiv -1\pmod 4
\end{cases}
\end{align*}
and
\begin{align*}
\vec{\tilde B}_t=\vec{B}_t+\sum_{i=0}^{t-1} \vec e_i=
\begin{cases}
\sum_{i=0}^{\frac{t-3}{2}} \vec e_{1+2i}, & \mbox{if } t\equiv  1\pmod 4 \\
\sum_{i=0}^{\frac{t-1}{2}} \vec e_{2i},   & \mbox{if } t\equiv -1\pmod 4 
\end{cases}
\end{align*}
Here $\vec e_i$'s are standard unit vectors in $\mathbb Z^t$ defined as 
$e_0=(1,0,\ldots,0),\ldots,\vec e_{t-1}=(0,\ldots,0,1)$.

We conclude this 
section with the following important 
observation.
If odd $t>1$, $k=0,1,\ldots,\frac{t-1}{2}$ and $\vec n\in\mathbb Z^t$, 
$\vec n\cdot\vec 1_t=0$, then
\beq
bg(\vec n)
=(-1)^{\frac{t-1}{2}}\left(\bigg\lfloor\frac{t}{4}\bigg\rfloor-k\right)
\label{eq:obs}
\eeq
iff $\vec n\equiv\vec B_t+\vec e_{i_0}+\vec e_{i_1}+\cdots+\vec e_{i_{2k}}
\pmod{2}$
for some $0\leq i_0 < i_1 <i_2<\cdots<i_{2k}\leq t-1$. 
In particular, if $\vec n\in\mathbb Z^t$, $\vec n\cdot\vec 1_t=0$, then
\beq
bg(\vec n)=(-1)^{\frac{t+1}{2}}\bigg\lfloor\frac{t+1}{4}\bigg\rfloor
\label{eq:obsa}
\eeq
iff $\vec n\equiv \vec{\tilde B}_t \pmod{2}$.
We leave the proof as an exercise for the interested reader.

\section{Combinatorial proof of $p_j(5n+4)\equiv 0\pmod 5$} \label{sec:4}

Throughout this section we assume that
\beqs
\vert\pi\vert\equiv 4\pmod 5
\eeqs
and
\beqs
\vert\pi_{\mbox{$5$-core}}\vert\equiv 4\pmod 5.
\eeqs
To prove (\ref{eq:*}) we shall require a few definitions. Following \cite{GKS}, we define the $\fc$ crank as
\beq
c_5(\pi):=2(r_0(\pi,5)-r_4(\pi,5))+(r_1(\pi,5)-r_3(\pi,5))+1\pmod 5.
\mylabel{eq:1*} 
\eeq
Note that if $\vert\pifc\vert\equiv4\pmod 5$, then obviously
\begin{align}
&n_0+n_1+n_2+n_3+n_4=0, \mylabel{eq:2*} \\ 
&n_1+2n_2+3n_3+4n_4\equiv 4\pmod 5.\mylabel{eq:3*} 
\end{align}
Here, $\vec n=(n_0,n_1,n_2,n_3,n_4)=\phi_2(\pifc)$.
Let's introduce a new vector $\vec \alpha(\vec n)=(\alpha_0,\alpha_1,\alpha_2,\alpha_3,\alpha_4)$, 
defined as
\begin{align}
&\alpha_0=\frac{n_0-3n_1-2n_2-n_3+1}{5}, \mylabel{eq:4*} \\ 
&\alpha_1=\frac{-3n_0-n_1-4n_2-2n_3+2}{5}, \mylabel{eq:5*} \\ 
&\alpha_2=\frac{-3n_0-n_1+n_2-2n_3+2}{5}, \mylabel{eq:6*} \\ 
&\alpha_3=\frac{n_0+2n_1+3n_2+4n_3+1}{5}, \mylabel{eq:7*} \\ 
&\alpha_4=\frac{4n_0+3n_1+2n_2+n_3-1}{5}. \mylabel{eq:8*} 
\end{align}
Using (\ref{eq:2*}, \ref{eq:3*}) it is easy to verify that $\vec\alpha(\vec n)\in\mathbb Z^5$
and that 
\beq
(\alpha_0+\alpha_1+\alpha_2+\alpha_3+\alpha_4)=1. 
\mylabel{eq:9*} 
\eeq
Inverting (\ref{eq:4*}--\ref{eq:8*}) we find that
\begin{align}
&n_0=\alpha_0+\alpha_4, \mylabel{eq:10*} \\ 
&n_1=-\alpha_0+\alpha_1+\alpha_4, \mylabel{eq:11*} \\ 
&n_2=-\alpha_1+\alpha_2, \mylabel{eq:12*} \\ 
&n_3=-\alpha_2+\alpha_3-\alpha_4, \mylabel{eq:13*} \\ 
&n_4=-\alpha_3-\alpha_4, \mylabel{eq:14*} 
\end{align}
Note that in terms of these new variables we have
\beq
c_5(\pi)\equiv\sum_{i=0}^4i\alpha_i\pmod 5,
\mylabel{eq:15*} 
\eeq
\beq
\vert\pi\vert=5Q(\vec\alpha)-1+5\sum_{i=0}^4\vert\pihat_i\vert,
\mylabel{eq:16*} 
\eeq
and
\begin{align}
\mbgr(\pi)&=\frac{1-(-1)^{\alpha_0+\alpha_1}-(-1)^{\alpha_1+\alpha_2}-\cdots-(-1)^{\alpha_4+\alpha_0}}{4} \nonumber \\
&+(-1)^{\alpha_0+\alpha_4}\mbgr(\pihat_0) \nonumber \\
&+(-1)^{\alpha_2+\alpha_3}\mbgr(\pihat_1) \nonumber \\
&+(-1)^{\alpha_1+\alpha_2}\mbgr(\pihat_2) \nonumber \\
&+(-1)^{\alpha_0+\alpha_1}\mbgr(\pihat_3) \nonumber \\
&+(-1)^{\alpha_3+\alpha_4}\mbgr(\pihat_4).
\mylabel{eq:17*} 
\end{align}
Here $\phi_1(\pi)=(\pifc,\pihat_0,\ldots,\pihat_4)$ and $Q(\vec\alpha):=\vec\alpha\cdot\vec\alpha-(\alpha_0\alpha_1+\alpha_1\alpha_2+\cdots+\alpha_4\alpha_0)$.
It is convenient to combine $\phi_1,\phi_2,\vec\alpha$ into a new invertible function $\Phi$, defined as
\beqs
\Phi(\pi)=(\vec\alpha(\phi_2(\pifc)),\vec\pihat),
\eeqs
where $\vec\pihat:=(\pihat_0,\ldots,\pihat_4)$. Following \cite{BG} we define
\begin{align*}
\widehat C_1(\vec\alpha)&=(\alpha_4,\alpha_0,\alpha_1,\alpha_2,\alpha_3), \\
\widehat C_2(\vec\pihat)&=(\pihat_4,\pihat_2,\pihat_3,\pihat_0,\pihat_1), \\
\widehat O(\pi)&=\Phi^{-1}(\widehat C_1(\vec\alpha),\widehat C_2(\vec\pihat)).
\end{align*}
We observe that operator $\widehat O$ has the following 
properties
\begin{align}
\big\vert\widehat O(\pi)\big\vert &=\vert\pi\vert, \nonumber\\
\widehat O^5(\pi)&=\pi, \nonumber\\
\mbgr\left(\widehat O(\pi)\right)&=\mbgr(\pi), \nonumber\\
c_5\left(\widehat O(\pi)\right)&\equiv 1+c_5(\pi)\pmod 5.
\mylabel{eq:18*} 
\end{align}
Clearly, $\widehat O$ preserves the norm and the $\mbgr$ of the partition. And so we can assemble all partitions of
$5n+4$ with $\mbgr=j$ into disjoint orbits:
\beqs
\pi, \quad \widehat O(\pi), \quad \widehat O^2(\pi), \quad \widehat O^3(\pi), \quad \widehat O^4(\pi).
\eeqs
Here, $\pi$ is some partition of $5n+4$ with $\mbgr=j$. Formula (\ref{eq:18*}) suggests that all five members of the same orbit are distinct. Clearly,
\beqs
p_j(5n+4)=5 \cdot\mbox{(number of orbits)}.
\eeqs
Hence, $p_j(5n+4)\equiv 0\pmod 5$, as desired.
In fact, we have the following
\begin{theorem}
\label{theoremfcc}
Let $j$ be any fixed integer. The residue of the $5$-core crank mod $5$
divides the partitions enumerated by $p_j(5n+4)$ into five equal classes.
\end{theorem}
We note that this theorem generalizes Theorem 4.1 \cite[p.717]{BG}.

\section{Identities for odd $t$-cores with extreme \mbgr \ values} \label{sec:5}

The main object of this section is to provide a proof of formulas 
(\ref{eq:**}) and (\ref{eq:***}).
Thoughout this section $t$ is presumed to be a positive odd integer. We will prove (\ref{eq:***}) first.
To this end we employ 
the observation \eqn{obs}
together with (\ref{eq:R}) to rewrite it as
\beq
\sum_{\substack{\vec n\in\mathbb Z^t,\vec n\cdot\vec 1_t=0 \\ \vec n\equiv\vec{\tilde B}_t\pmod{2}}}q^{\tilde Q(\vec n)}=q^{\frac{t^2-1}{8}}\frac{E^t(q^{4t})}{E(q^4)},
\mylabel{eq:1**} 
\eeq
where 
\beq
\tilde Q(\vec n):=\frac{t}{2}\vec n\cdot\vec n+\vec b_t\cdot\vec n.
\mylabel{eq:2**} 
\eeq
Next we introduce new summation variables $\vec{\tilde n}=(\tilde n_0,\ldots,\tilde n_{t-1})\in\mathbb Z^t$
as follows
\beq
\vec n=2\vec{\tilde n}+\sum_{i=0}^{\lfloor\frac{t-3}{4}\rfloor}
\left(\vec e_{\frac{t-3}{2}-2i}-\vec e_{\frac{t+1}{2}+2i}\right).
\mylabel{eq:3**} 
\eeq
Obviously, $\vec{\tilde n}$ is subject to the constraint
\beq
\vec{\tilde n}\cdot\vec 1_t=0.
\mylabel{eq:4**} 
\eeq
Note that in terms of new variables we have
\beq
\tilde Q(\vec n)=\tilde Q(\vec n)+(t-1)\vec 1_t\cdot\vec{\tilde n}=\frac{t^2-1}{8}
+4\big\{\frac{t}{2}\vec{\tilde n}\cdot\vec{\tilde n}+
\sigma_1+\sigma_2+\sigma_3\big\},
\mylabel{eq:5**} 
\eeq
where
\begin{align*}
\sigma_1 &= \sum_{i=0}^{\lfloor\frac{t-3}{4}\rfloor}(t-1-i)\tilde n_{\frac{t-3}{2}-2i}, \\
\sigma_2 &= \sum_{i=0}^{\lfloor\frac{t-3}{4}\rfloor}i\tilde n_{2i+\frac{t+1}{2}}, \\
\sigma_3 &= \sum_{i=-\lfloor\frac{t-1}{4}\rfloor}^{\lfloor\frac{t-1}{4}\rfloor}(\frac{t-1}{2}+i)
\tilde n_{\frac{t-1}{2}+2i}.
\end{align*}
At this point it is natural to perform further changes:
\begin{align*}
\tilde n_{\frac{t-3}{2}-2i} & \rightarrow\tilde n_{t-1-i}, \quad & 0\leq i\leq\bigg\lfloor\frac{t-3}{4}\bigg\rfloor \\
\tilde n_{\frac{t+1}{2}+2i} & \rightarrow\tilde n_i, \quad & 0\leq i\leq\bigg\lfloor\frac{t-3}{4}\bigg\rfloor \\
\tilde n_{\frac{t-1}{2}+2i} & \rightarrow\tilde n_{\frac{t-1}{2}+i}, \quad & 
-\bigg\lfloor\frac{t-1}{4}\bigg\rfloor\leq i\leq\bigg\lfloor\frac{t-1}{4}\bigg\rfloor.
\end{align*}
This way we obtain 
\begin{align*}
&\tilde Q(\vec n)=\frac{t^2-1}{8} + 4\tilde Q(\vec{\tilde n}), \\
&\vec{\tilde n}\in\mathbb Z^t, \quad \vec{\tilde n}\cdot\vec 1_t=0.
\end{align*}
And so with the aid of the Klyachko identity (\ref{eq:++}) we find that
\beq
C_{t,(-1)^{\frac{t+1}{4}}\lfloor\frac{t+1}{4}\rfloor} (q)=
\sum_{\substack{\vec{\tilde n}\in\mathbb Z^t \\ \vec{\tilde n}\cdot\vec 1_t=0}}
q^{\frac{t^2-1}{8}+4\tilde Q(\vec{\tilde n})}=q^{\frac{t^2-1}{8}}\frac{E^t(q^{4t})}{E(q^4)},
\mylabel{eq:6**} 
\eeq
as desired. To prove (\ref{eq:**}) we shall require the following lemma.
\begin{lemma}
\label{lemma1}
For a positive odd $t$
\beq
\psi^2(q^2)=q^{\frac{t-1}{2}}\psi^2(q^{2t})+\frac{E^3(q^{4t})}{f(-q^t,-q^{3t})}\sum_{i=0}^{\frac{t-3}{2}}q^i
\frac{f(q^{t-1-2i},-q^{1+2i})} {f(-q^{4i+2},-q^{4t-2-4i})}
\mylabel{eq:7**} 
\eeq
holds.
\end{lemma}

In the above we employed the Ramanujan notations
\begin{align}
\psi(q)&:=\frac{E^2(q^2)}{E(q)}=\sum_{n\ge 0}q^{\binom{n+1}{2}}, \mylabel{eq:8**} \\ 
f(a,b)&:=(ab,-a,-b;ab)_\infty. \mylabel{eq:9**} 
\end{align}
Using (\ref{eq:J}) we can easily show that 
\beq
f(a,b)=\sum_{n=-\infty}^\infty a^{\frac{n(n+1)}{2}} b^{\frac{n(n-1)}{2}}.
\mylabel{eq:10**} 
\eeq
Setting $a=q^{t-1-2i}, b=-q^{1+2i}, 0\leq i\leq\frac{t-3}{2}$ in (\ref{eq:10**}) and
dissecting we obtain
\begin{align}
f(q^{t-1-2i},-q^{1+2i})&=f(-q^{2+t+4i},-q^{3t-2-4i}) \nonumber \\
&+q^{t-1-2i}f(-q^{2-t+4i},-q^{5t-2-4i}).
\mylabel{eq:11**} 
\end{align}

To prove the above lemma we start with the Ramanujan $\mbox{}_1\psi_1$-summation formula \cite[II.29]{GR}
\beq
\sum_{n=-\infty}^\infty\frac{(a)_n}{(b)_n}z^n=
\frac{(az,\frac{q}{az},q,\frac{b}{a};q)_\infty}{(z,\frac{b}{az},b,\frac{q}{a};q)_\infty},
\;\;\;\;\;|\frac{b}{a}|<|z|<1.
\mylabel{eq:12**} 
\eeq
We set $b=aq$ to obtain
\beq
\sum_{n=-\infty}^\infty\frac{z^n}{1-aq^n}=
\frac{(az,\frac{q}{az},q,q;q)_\infty}{(z,\frac{q}{z},a,\frac{q}{a};q)_\infty}=
\frac{E^3(q)f(-az,-\frac{q}{az})}{f(-z,-\frac{q}{z})f(-a,-\frac{q}{a})},\;\;\;|q|<|z|<1.
\mylabel{eq:13**} 
\eeq
If we replace $q\rightarrow q^4,z=q,a=q^2$ in (\ref{eq:13**}) we find that
\beq
\sum_{n=-\infty}^\infty\frac{q^n}{1-q^{2+4n}}=\psi^2(q^2).
\mylabel{eq:14**} 
\eeq
Next we split the sum on the left of (\ref{eq:14**}) as
\beq
\psi^2(q^2)=\sum_{\substack{i=0 \\ i\neq\frac{t-1}{2}}}^{t-1} \sum_{m_i=-\infty}^\infty
q^i\frac{q^{tm_i}}{1-q^{2+4i}q^{4tm_i}}+
\sum_{m=-\infty}^\infty q^{\frac{t-1}{2}}\frac{q^{tm}}{1-q^{2t}q^{4tm}}.
\mylabel{eq:15**} 
\eeq
Using (\ref{eq:14**}) with $q\rightarrow q^t$ it is easy to recognize the last sum in
(\ref{eq:15**}) as $q^{\frac{t-1}{2}}\psi^2(q^{2t})$. And so we have
\beq
\psi^2(q^2)=q^{\frac{t-1}{2}}\psi^2(q^{2t})+\frac{E^3(q^{4t})}{f(-q^t,-q^{3t})}
\sum_{\substack{i=0 \\ i\neq\frac{t-1}{2}}}^{t-1}q^i
\frac{f(-q^{2+4i+t},-q^{3t-2-4i})}{f(-q^{2+4i},-q^{4t-2-4i})},
\mylabel{eq:16**} 
\eeq
where we have made a multiple use of (\ref{eq:13**}). 
Finally, folding the last sum in half and using (\ref{eq:11**}) we arrive at
\begin{align}
\psi^2(q^2)&=q^{\frac{t-1}{2}}\psi^2(q^{2t})+\sum_{i=0}^{\frac{t-3}{2}}
\frac{E^3(q^{4t})q^i}{f(-q^t,-q^{3t})f(-q^{2+4i},-q^{4t-2-2i})} \nonumber \\
&\times\big\{f\left(-q^{2+4i+t},-q^{3t-2-4i}\right)+q^{t-1-2i}f\left(-q^{5t-2-4i},-q^{2-t+4i}\right)\big\} \nonumber \\
&=q^{\frac{t-1}{2}}\psi^2(q^{2t})+\frac{E^3(q^{4t})}{f(-q^t,-q^{3t})}
\sum_{i=0}^{\frac{t-3}{2}}q^i\frac{f(q^{t-1-2i},-q^{1+2i})}{f(-q^{2+4i},-q^{4t-2-4i})}.
\mylabel{eq:17**} 
\end{align}
This concludes the proof of Lemma \ref{lemma1}.

We now move on to prove (\ref{eq:**}). Again, using the observation \eqn{obs},
we can rewrite it as 
\beq
\sum_{j=0}^{t-1}
\sum_{\substack{\vec n\in\mathbb Z^t,\vec n\cdot\vec 1_t=0 \\ \vec n\equiv\vec B_t+\vec e_j\pmod{2}}}
q^{\tilde Q(\vec n)}=q^{\frac{(t-1)(t-3)}{8}}F(t,q^2).
\mylabel{eq:18**} 
\eeq
Remarkably, (\ref{eq:18**}) is just the 
constant term in $z$ of the following more general identity
\begin{align}
&\sum_{j=0}^{t-1}
\sum_{\substack{\vec n\in\mathbb Z^t \\ \vec n\equiv\vec B_t+\vec e_j\pmod{2}}}
q^{\tilde Q(\vec n)}z^{\frac{\vec n\cdot\vec 1_t}{2}} \nonumber \\
&=q^{\frac{(t-1)(t-3)}{8}}F(t,q^2)
\sum_{n=-\infty}^\infty q^{2n^2+(t-1)n}z^n.
\mylabel{eq:19**} 
\end{align}
To prove (\ref{eq:19**}) we observe that its right hand side satisfies the {\it first} order functional equation
\beq
\widehat D_{t,q}(f(z))=f(z),
\mylabel{eq:20**} 
\eeq
where
\beqs
\widehat D_{t,q}(f(z)):=zq^{t+1}f(zq^4).
\eeqs
After a bit of labor one can verify that for $0\leq i\leq t-1$
\beq
\widehat D_{t,q}
\left(\sum_{\substack{\vec n\in\mathbb Z^t \\ \vec n\equiv\vec B_t+\vec e_i\pmod{2}}}
q^{\tilde Q(\vec n)}z^{\frac{\vec n\cdot\vec 1_t}{2}}\right)
=\sum_{\substack{\vec n\in\mathbb Z^t \\ \vec n\equiv\vec B_t+\vec e_{i+2}\pmod{2}}}
q^{\tilde Q(\vec n)}z^{\frac{\vec n\cdot\vec 1_t}{2}},
\mylabel{eq:21**} 
\eeq
where $\vec e_t:=\vec e_0$ and $\vec e_{t+1}:=\vec e_1$. Clearly, (\ref{eq:21**}) implies that 
the left hand side of (\ref{eq:19**}) satisfies (\ref{eq:20**}), as well.
It remains to verify (\ref{eq:19**}) at one nontrivial point. To this end we set
\begin{align*}
z=
\begin{cases}
1,   & \mbox{if } t\equiv-1\pmod 4, \\
q^2, & \mbox{if } t\equiv 1\pmod 4
\end{cases}
\end{align*}
in (\ref{eq:19**}), and then replace $q^2\rightarrow q$ to get with the help of (\ref{eq:J})
\begin{align}
&q^{\frac{t-1}{2}}\psi(q^{2t})
\prod_{j=0}^{\frac{t-3}{2}}f^2\left(q^{1+2j},q^{2t-1-2j}\right) \nonumber \\
&\times\bigg\{1+\sum_{i=1}^{\frac{t-1}{2}}q^{-i}
\frac{f(q^t,q^t)f(q^{2i},q^{2t-2i})}{\psi(q^{2t})f(q^{t+2i},q^{t-2i})}\bigg\}
=\psi(q^2)F(t,q).
\mylabel{eq:22**} 
\end{align}
To proceed further we need to verify two product identities
\beqs
\psi(q^2)\prod_{j=0}^{\frac{t-3}{2}}f^2\left(q^{1+2j},q^{2t-1-2j}\right)=\psi(q^{2t})F(t,q)
\eeqs
and
\beqs
\psi(q^{2t})\frac{f(q^t,q^t)f(q^{2i},q^{2t-2i})}{f(q^{t+2i},q^{t-2i})}=
E^3(q^{4t})\frac{f(q^{2i},-q^{t-2i})}{f(-q^t,-q^{3t})f(-q^{2t+4i},-q^{2t-4i})}, \quad i\in\mathbb N.
\eeqs
Next, we multiply both sides of (\ref{eq:22**}) by $\frac{\psi(q^2)}{F(t,q)}$
and simplify to arrive at
\beq
q^{\frac{t-1}{2}}\psi^2(q^{2t})+\frac{E^3(q^{4t})}{f(-q^t,-q^{3t})}
\sum_{i=1}^{\frac{t-1}{2}}q^{\frac{t-1}{2}-i}
\frac{f(q^{2i},-q^{t-2i})}{f(-q^{2t+4i},-q^{2t-4i})}=\psi^2(q^2),
\mylabel{eq:23**} 
\eeq
which is essentially the identity in Lemma \ref{lemma1}. This concludes our proof of (\ref{eq:19**}).
It follows that (\ref{eq:18**}), (\ref{eq:**}) hold true.

\section{$5$-cores with prescribed  \mbgr} \label{sec:6}

Formula (\ref{eq:bound}) suggests that $\mbgr(\pifc)$ can assume just three values: $0,\pm1$. This means that
\beq
a_5(n)=a_{5,-1}(n)+a_{5,0}(n)+a_{5,1}(n).
\mylabel{eq:*0} 
\eeq
The generating function of version (\ref{eq:*0}) is
\beq
\frac{E^5(q^5)}{E(q)}=C_{5,-1}(q)+C_{5,0}(q)+C_{5,1}(q).
\mylabel{eq:*1} 
\eeq
In the last section we proved (\ref{eq:**}), (\ref{eq:***}). These identities with $t=5$ state that
\begin{align}
C_{5,-1}(q)&=q^3\frac{E^5(q^{20})}{E(q^4)}, \mylabel{eq:*2} \\ 
C_{5,1}(q) &=qF(5,q^2). \mylabel{eq:*3} 
\end{align}
By (\ref{eq:aform}) we observe that $C_{t,j}(q)$ is either an odd
or an even functions of $q$ with parity determined by the parity
of $j$. Therefore, $C_{5,0}(q)$ is an even function of $q$,
and $C_{5,\pm1}(q)$ are odd functions of $q$.
Consequently, we see that
\beq
\epart\left(\frac{E^5(q^5)}{E(q)}\right)=C_{5,0}(q)
\mylabel{eq:*5} 
\eeq
where
\beqs
\epart(f(x)):=\frac{f(x)+f(-x)}{2}.
\eeqs
In this section we will show that $C_{5,0}(q)$ can be expressed as a sum of two infinite products
\beq
C_{5,0}(q)=R(q^2),
\mylabel{eq:*6} 
\eeq
where
\beq
R(q):=\frac{E^4(q^{10})E(q^5)E^2(q^4)}{E^2(q^{20})E(q)}
+q\frac{E^2(q^{20})E^3(q^5)E^6(q^2)}{E^2(q^{10})E^2(q^4)E^3(q)}.
\mylabel{eq:*7} 
\eeq
It is easy to rewrite (\ref{eq:*7}) in a manifestly positive way as
\beqs
R(q)=f(q,q^4)f(q^2,q^3)\big\{\varphi(q^5)\psi(q^2)+q\varphi(q)\psi(q^{10})\big\},
\eeqs
where
\beqs
\varphi(q) := f(q,q) = \sum_{n=-\infty}^{\infty} q^{n^2} =
\frac{E^5(q^2)}{E^2(q^4)E^2(q)},
\eeqs
and $\psi(q)$ is defined in (\ref{eq:8**}).
Formula (\ref{eq:*6}) enabled us to discover and prove the new Lambert series identity
\beq
R(q)=\sum_{i=0}^1 \sum_{n=-\infty}^\infty (-1)^iq^{5n+i}
\frac{1+q^{1+2i+10n}}{(1-q^{1+2i+10n})^2}.
\mylabel{eq:*8} 
\eeq
In what follows we will require three identities:
\beq
\big[ux,\frac{u}{x},vy,\frac{v}{y};q\big]_\infty=
\big[uy,\frac{u}{y},vx,\frac{v}{x};q\big]_\infty+\frac{v}{x}
\big[xy,\frac{x}{y},uv,\frac{u}{v};q\big]_\infty,
\mylabel{eq:*9} 
\eeq
(\cite[ex. 5.21]{GR})
\beq
f(a,b)f(c,d)=f(ac,bd)f(ad,bc)+af\left(\frac{b}{c},\frac{c}{b}(abcd)\right)f\left(\frac{b}{d},\frac{d}{b}(abcd)\right),
\mylabel{eq:*10} 
\eeq
provided $ab=cd$ (\cite{BB}) and
\beq
\frac{E^5(q^5)}{E(q)}=\sum_{i=1}^2 \sum_{n=-\infty}^\infty (-1)^{i+1}
\frac{q^{5n+i-1}}{(1-q^{5n+i})^2},
\mylabel{eq:*11} 
\eeq
(\cite[ex. 5.7]{GR}, \cite[p.8]{GKS}).
Here
\begin{align*}
[a;q]_\infty &=\left(a,\frac{q}{a};q\right)_\infty, \\
[a_1,a_2,\ldots,a_n;q]_\infty &=\prod_{i=1}^n[a_i;q]_\infty.
\end{align*}
Next, we wish to establish the validity of
\beq
F(5,q)=\frac{E(q^{10})E^2(q^5)E^3(q^2)}{E^2(q)}
=\frac{E^5(q^5)}{E(q)}+q\frac{E^5(q^{10})}{E(q^2)}.
\mylabel{eq:*12} 
\eeq
To this end we multiply both sides of (\ref{eq:*12}) by
\beqs
\frac{\big[q,q^3;q^{10}\big]_\infty^2\big[q^2,q^4;q^{10}\big]_\infty}{E^4(q^{10})}
\eeqs
to obtain after simplification that
\beq
[q^2,q^2,q^4,q^6;q^{10}]_\infty=[q,q^3,q^5,q^5;q^{10}]_\infty+q[q,q,q^3,q^3;q^{10}]_\infty.
\mylabel{eq:*13} 
\eeq
But the last equation is nothing else but (\ref{eq:*9}) with 
$q$ replaced by $q^{10}$ and
$u=q^2,v=q^5,x=1,y=q$.
We now combine
\beqs
\epart\left(q\frac{E^5(q^5)}{E(q)}\right)=qC_{5,-1}(q)+qC_{5,1}(q),
\eeqs
with (\ref{eq:*2}), (\ref{eq:*5}), and (\ref{eq:*12}) to obtain
\beq
\epart\left(q\frac{E^5(q^5)}{E(q)}\right)= 2q^4 \frac{E^5(q^{20})}{E(q^4)} + q^2\frac{E^5(q^{10})}{E(q^2)}.
\mylabel{eq:*14} 
\eeq
This can be stated as the following eigenvalue problem
\beq
T_{2}\left(q\frac{E^5(q^5)}{E(q)}\right)= q\frac{E^5(q^5)}{E(q)},
\mylabel{eq:*15} 
\eeq
where for prime $p$ the Hecke operator $T_{p}$ is defined by its action as
\beqs
T_{p}\left(\sum_{n\ge 0}a_nq^n\right)=\sum_{n\ge 0}a_{pn}q^n+p\legendre{p}{5}\sum_{n\ge 0}a_nq^{pn},
\eeqs
with $\legendre{a}{b}$ being the Legendre symbol. We remark that (\ref{eq:*15}) 
is the $p=2$ 
case of the more general formula 
\beq
T_{p}\left(q\frac{E^5(q^5)}{E(q)}\right)= \left(p+\legendre{p}{5}\right)\bigg(q\frac{E^5(q^5)}{E(q)}\bigg),
\mylabel{eq:*16} 
\eeq
which can be deduced from (\ref{eq:*11}). We shall not supply the details. Instead, we note that (\ref{eq:*16})
together with (\ref{eq:*2}, \ref{eq:*3}, \ref{eq:*5}) implies that
\beq
T_{\tilde p}(qC_{5,j}(q))=\left(\tilde p+\legendre{\tilde p}{5}\right)(qC_{5,j}(q)), \quad j=0\pm 1.
\mylabel{eq:*16a} 
\eeq
Here, $\tilde p$ is an odd prime.

To prove (\ref{eq:*6}) we use (\ref{eq:*12}) to deduce that 
\beq
\epart\left(\frac{E^5(q^5)}{E(q)}\right)= \epart(F(5,q))=E(q^{10})E^3(q^2)\cdot \epart\left(\frac{E^2(q^5)}{E^2(q)}\right).
\mylabel{eq:*17} 
\eeq
To proceed further we employ (\ref{eq:*10}) with $a=q,b=q^9,c=q^3,d=q^7$ to get
\begin{align}
\frac{E(q^5)}{E(q)} &= \frac{E(q^4)}{E(q^{20})E^2(q^2)} f(q,q^9)f(q^3,q^7) \nonumber \\
&=\frac{E(q^4)}{E(q^{20})E^2(q^2)}\{f(q^4,q^{16})f(q^8,q^{12})+q f(q^6,q^{14})f(q^2,q^{18})\}\nonumber \\
&=\frac{E^2(q^{20})E(q^8)}{E(q^{40})E^2(q^2)} + q\frac{E(q^{40})E(q^{10})E^3(q^4)}{E(q^{20})E(q^8)E^3(q^2)}.
\mylabel{eq:*18} 
\end{align}
It is clear that
\beq
\epart\left(\frac{E^2(q^5)}{E^2(q)}\right)=\frac{E^4(q^{20})E^2(q^8)}{E^2(q^{40})E^4(q^2)}
+q^2\frac{E^2(q^{40})E^2(q^{10})E^6(q^4)}{E^2(q^{20})E^2(q^8)E^6(q^2)}.
\mylabel{eq:*19} 
\eeq
Combining (\ref{eq:*17}) and (\ref{eq:*19}) we find that
\beq
\epart\left(\frac{E^5(q^5)}{E(q)}\right)=R(q^2).
\mylabel{eq:*20} 
\eeq
The last formula together with (\ref{eq:*5}) implies (\ref{eq:*6}). 
Next, we rewrite (\ref{eq:*11}) as
\beqs
\frac{E^5(q^5)}{E(q)}=\sum_{i=1}^2\sum_{n=-\infty}^\infty(-1)^{i+1}
\frac{q^{5n+i-1}(1+2q^{5n+i}+q^{10n+2i})}{(1-q^{10n+2i})^2}.
\eeqs
Clearly,
\begin{align}
\epart\left(\frac{E^5(q^5)}{E(q)}\right)&=\sum_{i=1}^2\sum_{\substack{n=-\infty \\ n\equiv i-1\pmod2}}^\infty
(-1)^{i+1}\frac{q^{5n+i-1}(1+q^{10n+2i})}{(1-q^{10n+2i})^2} \nonumber \\
&=\sum_{i=0}^1\sum_{n=-\infty}^\infty(-1)^i
\frac{q^{10n+i}(1+q^{20n+4i+2})}{(1-q^{20n+4i+2})^2}.
\mylabel{eq:*21} 
\end{align}
Formula (\ref{eq:*8}) with $q\rightarrow q^2$ follows easily from (\ref{eq:*20}) and (\ref{eq:*21}).
Before we move on we wish to summarize some of the above observations in the formula below
\begin{align}
\frac{E^5(q^5)}{E(q)} &= \Big\{\frac{E^4(q^{20})E(q^{10})E^2(q^8)}{E^2(q^{40})E(q^2)}
+q^2\frac{E^2(q^{40})E^3(q^{10})E^6(q^4)}{E^2(q^{20})E^2(q^8)E^3(q^2)}\Big\}\nonumber \\
&+q\Big\{\frac{E^5(q^{10})}{E(q^2)}+2q^2\frac{E^5(q^{20})}{E(q^4)}\Big\}.
\mylabel{eq:*21a} 
\end{align}
In \cite{GKS}, the authors used (\ref{eq:*11}) to find explicit formulas 
for the coefficients
\beq
a_5(n)=\frac{2^{d+1}+(-1)^d}{3}\cdot 5^c\cdot\prod_{i=1}^s\frac{p_i^{a_i+1}-1}{p_i-1}
\prod_{j=1}^t\frac{q_j^{b_j+1}+(-1)^{b_j}}{q_j+1}.
\mylabel{eq:*22} 
\eeq
Here
\beq
n+1=2^d5^c\prod_{i=1}^s p_i^{a_i}\prod_{j=1}^t q_j^{b_j}
\mylabel{eq:*23} 
\eeq
is the prime factorization of $n+1$ and $p_i\equiv\pm1\pmod5, 1\le i\le s$ and $q_j\equiv\pm2\pmod5, 1\le j\le t$
are odd primes. Formulas (\ref{eq:*2}), (\ref{eq:*3}), (\ref{eq:*5}) and (\ref{eq:*12}) suggest the following relations.
For $n\in\mathbb N$ and $r=0,1,2,3$ one has
\begin{align}
a_{5,0}(n)&=
\begin{cases}
a_5(n), &\mbox{if } n\equiv0\pmod2, \\
0, &\mbox{otherwise},
\mylabel{eq:*24} 
\end{cases} \\
a_{5,-1}(4n+r)&=
\begin{cases}
a_5(n), &\mbox{if } r=3, \\
0, &\mbox{otherwise}, 
\mylabel{eq:*25} 
\end{cases} \\
a_{5,1}(4n+r)&=
\begin{cases}
a_5(2n), &\mbox{if } r=1, \\
a_5(n)+a_5(2n+1), &\mbox{if } r=3, \\
0, &\mbox{if } r=0,2. 
\mylabel{eq:*26} 
\end{cases}
\end{align}
These relations together with (\ref{eq:*22}) enabled us to derive explicit formulas for $a_{5,j}(n)$ with $-1\le j\le 1$.
In particular, if the prime factorization of $n+1$ is given by (\ref{eq:*23}), then
\beq
a_{5,1}(4n+3)=2^{d+1}5^c\prod_{i=1}^s\frac{p_i^{a_i+1}-1}{p_i-1}
\prod_{j=1}^t\frac{q_j^{b_j+1}+(-1)^{b_j}}{q_j+1}.
\mylabel{eq:*27} 
\eeq
We would like to conclude this section with the following discussion. It is easy to check that (\ref{eq:*16a})
implies that
\beq
a_{5,j}(pn+p-1)+p\legendre{p}{5}a_{5,j}\left(\frac{n+1}{p}-1\right)
=\left(p+\legendre{p}{5}\right)a_{5,j}(n), \quad j=0,\pm1,
\mylabel{eq:*28} 
\eeq
where $p$ is odd prime, $n\in\mathbb N$ and $a_{5,j}(x)=0$ if $x\not\in\mathbb Z$. Setting $p=5$ we find that
\beq
a_{5,j}(5n+4)=5a_{5,j}(n), \quad j=0,\pm1.
\mylabel{eq:*29} 
\eeq
This is a refinement of the well-known result
\beq
a_5(5n+4)=5a_5(n),
\mylabel{eq:*30} 
\eeq
proven in \cite{GKS}. 
We can prove (\ref{eq:*29}) by adapting the combinatorial proof in \cite{GKS}.

Let's define
\beqs
\vec n=(n_0,n_1,n_2,n_3,n_4)=\phi_2(\pifc)
\eeqs
for some $\pifc$ with $\mbgr(\pifc)=j$ and $\vert\pifc\vert=n$. 
Consider map $\vec n\rightarrow\vec{\tilde n}=(\tilde n_0,\tilde n_1,\tilde n_2,\tilde n_3,\tilde n_4)$
with 
\begin{align*}
\tilde n_0 &= n_1+2n_2+2n_4+1, \\
\tilde n_1 &= -n_1-n_2+n_3+n_4+1, \\
\tilde n_2 &= 2n_1+n_2+2n_3, \\
\tilde n_3 &= -2n_2-2n_3-n_4-1, \\
\tilde n_4 &= -2n_1-n_3-2n_4-1. \\
\end{align*}
Obviously $\vec{\tilde n}\in\mathbb Z^5$ and $\vec{\tilde n}\cdot\vec 1_5=0$ and so we can define $\tilde\pi_\fc=\phi_2^{-1}(\vec{\tilde n})$. It is easy to check that
\beqs
\vert\tilde\pi_\fc\vert=5n+4,
\eeqs
and that
\beqs
\mbgr(\tilde\pi_\fc)=j,
\eeqs
and
\beqs
c_5(\tilde\pi_\fc)\equiv 4\pmod5.
\eeqs
Recall that the orbit $\{\tilde\pi_\fc,\widehat O(\tilde\pi_\fc),\ldots,\widehat O^4(\tilde\pi_\fc)\}$ contains just one member with $c_5\equiv4\pmod5$. 
And so each $5$-core of $n$ with $\mbgr$ $j$ is in $1-1$ correspondence with 
an appropriate $5$-member orbit of $t$-cores of $5n+4$ with $\mbgr$ $j$. 
This observation yields a combinatorial proof of (\ref{eq:*29}).

\section{Outlook}

Given our 
combinatorial proof of
\beqs
p_j(5n+4)\equiv 0\pmod5,\quad j\in\mathbb Z
\eeqs
one may wonder about a combinatorial proof of the other mod $5$ 
congruences (\ref{eq:10}-\ref{eq:40}). We strongly suspect that such proof will be dramatically different from the one discussed in Section 4. In addition, one would like to have combinatorial insights into (\ref{eq:*28}) for $p\neq 5$.

In this paper we found ``positive'' $eta$-quotient representations for $C_{5,j}(q),-1\leq j\leq1$. In the general case 
(odd  $t$, $-\lfloor\frac{t-1}{4}\rfloor\le j\le\lfloor\frac{t+1}{4}\rfloor)$, we established such representation only for $C_{t,\pm\lfloor\frac{t\pm1}{4}\rfloor}(q)$. Clearly, one wants to find ``positive'' $eta$-quotient representations for other admissible values of $\mbgr$. (See \cite{BY} for a 
fascinating discussion of the $t=7$ case).

Finally, we observe that (\ref{eq:bgd'}) is the $s=2$ case of the following more 
general definition
\beqs
\mbox{gbg-rank}(\pi,s)=\sum^{s-1}_{j=0}r_j(\pi,s)\omega_s^j
\eeqs
with
\beqs
\omega_s=e^{i\frac{2\pi}{s}}.
\eeqs

Many identities, proven here, can be generalized further. For example, 
we can prove that if
$(s,t)=1$ then
\beq
\mbox{gbg-rank}(\pi_{\tc},s)=
\frac{\sum_{i=0}^{t-1}\omega_s^{i+1}(\omega_s^{tn_i}-1)}{(1-\omega_s^t)(1-\omega_s)}
\mylabel{eq:*31} 
\eeq
and for $1\leq i\leq s-1$ that
\beq
\sum_{\mbox{gbg-rank}(\pi_{t\mbox{-core}},s)=g(i)}q^{|\pi_{t\mbox{-core}}|}=q^{a(i)}F_i(q^s).
\mylabel{eq:*32} 
\eeq
Here, 
\beqs
(n_0,n_1,\ldots,n_{t-1})=\phi_2(\pi_{t\mbox{-core}}),
\eeqs

\beqs
a(i)=\frac{(t^2-1)(s^2-1)}{24}-\frac{(t-1)(s-i)i}{2},
\eeqs
\beqs
g(i)=\frac{1}{(1-\omega_s)(1-\frac{1}{\omega_s})}-\omega_s^{\frac{t-1}{2}}\frac{1+\frac{t-1}{\omega_s^i}}{(1-\omega_s^t)(1-\frac{1}{\omega_s})},
\eeqs
\beqs
F_i(q)=E(q^s)E(q^{st})^{t-2}\frac{[q^{it};q^{st}]_\infty}{[q^i;q^s]_\infty}.
\eeqs
Setting $s=2$ in (\ref{eq:*31}), (\ref{eq:*32}) we obtain (\ref{eq:2}), (\ref{eq:**}), respectively.

In addition we can show that
\beq
\sum_{\mbox{gbg-rank}(\pi_{t\mbox{-core}},s)=g(0)}q^{|\pi_{t\mbox{-core}}|}=q^{a(0)}\frac{E(q^{s^2t})^t}{E(q^{s^2})}.
\mylabel{eq:*33} 
\eeq
Setting $s=2$ in (\ref{eq:*33}) we get (\ref{eq:***}). 

In \cite{OS} Olsson and Stanton defined so-called $(s,t)$-good partitions. Surprisingly, $t$-cores with 
$\mbox{gbg-rank}=g(0)$ coincide with $(t,s)$-good partitions.

Let $\nu(t,s)$ denote a number of distinct values that $\mbox{gbg-rank}(\pi_{t\mbox{-core}},s)$ may assume.
Then it can be shown that 
\beqs
\nu(s,t)\leq\frac{\binom{t+s}{t}}{t+s},
\eeqs
provided that $(s,t)=1$.
Morever, if $s$ is prime or if $s$ is a composite number and $t<2p$ then
\beqs
\nu(s,t)=\frac{\binom{t+s}{t}}{t+s}.
\eeqs
Here, $p$ is a smallest prime divisor of $s$ and $(s,t)=1$.

Details of these and related results will be left to a later paper.

\bigskip
\noindent
\textbf{Acknowledgement}

We would like to thank Robin Chapman, Ole Warnaar, Herbert Wilf and 
Hamza Yesilyurt for 
their kind interest and stimulating discussions. 

\bibliographystyle{amsplain}

\end{document}